\documentclass[11pt,a4paper]{article}

\usepackage{amsmath}
\usepackage{amssymb,verbatim}
\usepackage{amsfonts}
\usepackage{amsthm}
\usepackage{pifont}
\usepackage{graphicx}

\newcommand{\beq}{\begin{equation*}}
\newcommand{\eeq}{\end{equation*}}
\newcommand{\beqn}{\begin{equation}}
\newcommand{\eeqn}{\end{equation}}

\newcommand{\dd}{\mathrm{d}}
\newcommand{\ii}{\mathrm{i}}
\newcommand{\ee}{\mathrm{e}}
\newcommand{\sinc}{\mathrm{sinc}}

\newcommand{\ds}{\displaystyle}

\begin{document}

\title{A remark concerning sinc integrals}
\author{Uwe B\"asel}
\date{}
\maketitle

\begin{abstract}
\noindent We give a simple proof of Hanspeter Schmid's result that
\beq
  K_n:=2\int_0^\infty\cos t\,\prod_{k=0}^n\sinc\left(\frac{t}{2k+1}\right)\dd t=\frac{\pi}{2}
  \;\;\,\mbox{if}\;\;\,
  n\in\{0,1,\ldots,55\}\,,
\eeq
and $K_n<\pi/2$ if $n\geq 56$. Furthermore, we present two sinc integrals where the value $\pi/2$ is undercut as soon as $n\geq 418$ and $n\geq 3091$, respectively.\\[0.2cm]
\textbf{Mathematics Subject Classification:} 33B10, 26D15\\[0.2cm]
\textbf{Keywords:} sinc integrals, improper integrals
\end{abstract}

\section{Introduction}

We define
\beq
  \sinc(t) := \left\{\begin{array}{c@{\quad\mbox{if}\quad}c}
	\dfrac{\sin t}{t} & t\not=0\,,\\[0.3cm]
	1 & t=0\,.
  \end{array}\right.
\eeq
Recently, Schmid \cite{Schmid} investigated the integrals
\beqn \label{tau_n}
  \tau_n:=\int_0^\infty\prod_{k=0}^n\sinc(a_k t)\,\dd t
\eeqn
and
\beqn \label{eps_n}
  \varepsilon_n:=\int_0^\infty 2\cos a_0 t\,\prod_{k=0}^n\sinc(a_k t)\,\dd t
\eeqn
with positive real numbers $a_0\geq a_1\geq\ldots\geq a_n$, and proved that
\beqn \label{tau_na}
  \left.\begin{array}{lll}
	\tau_0 = \dfrac{\pi}{2a_0} \,, &
	\tau_n = \dfrac{\pi}{2a_0} & \mbox{for all $n$ with} \quad 
	\ds{\sum_{k=1}^n a_k\leq a_0}\,,\\[0.5cm]
	& \tau_n < \tau_{n-1} & \mbox{for all other $n$}\,, 
  \end{array}\right\}
\eeqn 
and
\beqn \label{eps_na}
  \left.\begin{array}{lll}
	\varepsilon_0 = \dfrac{\pi}{2a_0}\,, &
	\varepsilon_n = \dfrac{\pi}{2a_0} & \mbox{for all $n$ with} \quad 
	\ds{\sum_{k=1}^n a_k\leq 2a_0}\,,\\[0.5cm]
	& \varepsilon_n < \varepsilon_{n-1} & \mbox{for all other $n$}\,. 
  \end{array}\right\}
\eeqn 
Result \eqref{eps_na} was found by considering point-symmetric properties of the Fourier transforms
\beqn \label{FT}
  F_n(\omega):=\int_{-\infty}^\infty\,\ee^{-\ii\omega t}\prod_{k=0}^n\sinc(a_k t)\,\dd t
\eeqn
(see \cite[pp.\:14-15,\:17]{Schmid}).\\[0.2cm]  
A special case of the integral \eqref{tau_n} is
\beqn \label{J_n}
  J_n:=\int_0^\infty\prod_{k=0}^n\sinc\left(\frac{t}{2k+1}\right)\dd t\,.
\eeqn
Therefore, from \eqref{tau_na} one easily finds \cite[pp.\:11-12,\:16]{Schmid} 
\beqn \label{J_na}
  \left.\begin{array}{l@{\quad\mbox{if}\quad}l}
	J_n=\pi/2 & n\in\{0,1,\ldots,6\}\,,\\[0.1cm]
	J_n<J_{n-1} & n\geq 7\,.
  \end{array}\right\}
\eeqn
This result was first proven by Borwein \& Borwein \cite[p.\:86]{Borwein_Borwein}.\\[0.2cm] 
A special case of the integral \eqref{eps_n} is
\beqn \label{K_n}
  K_n := 2\int_0^\infty\cos t\,\prod_{k=0}^n\sinc\left(\frac{t}{2k+1}\right)\dd t\,.
\eeqn
From \eqref{eps_na} one finds \cite[pp.\:12,\:17]{Schmid} 
\beqn \label{K_na}
  \left.\begin{array}{l@{\quad\mbox{if}\quad}l}
	K_n=\pi/2 & n\in\{0,1,\ldots,55\}\,,\\[0.1cm]
	K_n<K_{n-1} & n\geq 56\,.
  \end{array}\right\}
\eeqn
Due to the properties \eqref{J_na} and \eqref{K_na}, Schmid called \eqref{J_n} and \eqref{K_n} as {\em curious integrals}.
 
\section{Alternative proof of \eqref{eps_na}} \label{proof}

At first we show that the result \eqref{eps_na} can be obtained in a very simple way without using the Fourier transforms \eqref{FT} and its symmetry properties.\\[0.2cm]
Using $2\sin x\cos x=\sin 2x$, we have
\begin{align*}
  \varepsilon_n
	= {} & \int_0^\infty 2\cos(a_0t)\,\prod_{k=0}^n\sinc(a_k t)\:\dd t\\
	= {} & 2\int_0^\infty\cos(a_0t)\;\sinc(a_0t)\,\prod_{k=1}^n\sinc(a_k t)\:\dd t\displaybreak[0]\\
	= {} & 2\int_0^\infty\cos(a_0t)\;\frac{\sin(a_0t)}{a_0t}\,\prod_{k=1}^n\sinc(a_k t)\:\dd t\displaybreak[0]\\
	= {} & 2\int_0^\infty\frac{\sin(2a_0t)}{2a_0t}\,\prod_{k=1}^n\sinc(a_k t)\:\dd t\displaybreak[0]\\
	= {} & 2\int_0^\infty\sinc(2a_0t)\,\prod_{k=1}^n\sinc(a_k t)\:\dd t\,.
\end{align*}
Substituting $b_0=2a_0$, and $b_k=a_k$ if $k\in\{1,2,\ldots,n\}$, we have
\beq
  \varepsilon_n
	= 2\int_0^\infty\prod_{k=0}^n\sinc(b_k t)\:\dd t\,.
\eeq
From \eqref{tau_na} it follows that
\beqn \label{eps_nb}
 \left.\begin{array}{lll}
	\hspace*{-0.5cm}
	\varepsilon_0 = 2\cdot\dfrac{\pi}{2b_0} = \dfrac{\pi}{b_0}\,, &
	\varepsilon_n = 2\cdot\dfrac{\pi}{2b_0} = \dfrac{\pi}{b_0} & \mbox{for all $n$ with} \quad 
	\ds{\sum_{k=1}^n b_k\leq b_0}\,,\\[0.5cm]
	& \varepsilon_n < \varepsilon_{n-1} & \mbox{for all other $n$}\,. 
  \end{array}\right\}
\eeqn
Substituting back, \eqref{eps_na} follows immediately.\\[0.2cm] 
This shows that integral \eqref{eps_n} may be considered as special case of integral~\eqref{tau_n}, and result \eqref{eps_na} as special case of result \eqref{tau_na}.

\section{Proof of \eqref{K_na}}
  
From Section \ref{proof} we find that integral \eqref{K_n} may be written as
\beq
  K_n = 2\int_0^\infty\prod_{k=0}^n\sinc(b_kt)\,\dd t
\eeq
where
\beq
  b_0 := 2 \,,\quad b_k := \frac{1}{2k+1} \;\;\mbox{if}\;\; k\in\{1,2,\ldots,n\}\,. 
\eeq
Using \eqref{eps_nb} yields
\beq
 \left.\begin{array}{lll}
	K_0 = \dfrac{\pi}{2}\,, &
	K_n = \dfrac{\pi}{2} & \mbox{for all $n$ with} \quad 
	\ds{\sum_{k=1}^n\frac{1}{2k+1}\leq 2}\,,\\[0.5cm]
	& K_n < K_{n-1} & \mbox{for all other $n$}\,. 
  \end{array}\right.
\eeq
Result \eqref{K_na} follows.

\section{Many other curious integrals}

We put
\beq \label{I}
  I_n(b):=b\int_0^\infty\sinc(bt)\,\prod_{k=0}^n\sinc\left(\frac{t}{2k+1}\right)\dd t
\eeq
with a real number $b\geq 1$. From \eqref{tau_na} it is obviously clear that
\beq
  b\int_0^\infty\sinc(bt)\,\dd t = \frac{\pi}{2}\,,
\eeq
and
\beq \label{Ia}
  \left.\begin{array}{l@{\quad\mbox{if}\quad}l}
	I_n(b) = \dfrac{\pi}{2} &  
	\ds{\sum_{k=0}^n\frac{1}{2k+1}\leq b}\,,\\[0.5cm]
	I_n(b) < \dfrac{\pi}{2} & 
	\ds{\sum_{k=0}^n\frac{1}{2k+1}>b}\,. 
  \end{array}\right.
\eeq
Here are some results for special cases:
\begin{itemize}
\item $I_n(1)=\pi/2$ if $n=0$, $I_n(1)<\pi/2$ if $n\geq 1$, 
\item $I_n(2)=\pi/2$ if $n\in\{0,1,\ldots,6\}$, $I_n(2)<\pi/2$ if $n\geq 7$, 
\item $I_n(3)=\pi/2$ if $n\in\{0,1,\ldots,55\}$, $I_n(3)<\pi/2$ if $n\geq 56$, 
\item $I_n(4)=\pi/2$ if $n\in\{0,1,\ldots,417\}$, $I_n(4)<\pi/2$ if $n\geq 418$, 
\item $I_n(5)=\pi/2$ if $n\in\{0,1,\ldots,3090\}$, $I_n(5)<\pi/2$ if $n\geq 3091$.
\end{itemize}
A calculation with Mathematica gave
\begin{align*}
  I_7(2)
	= {} & \frac{168579263752211300739165075916829279}
			{337158527504429357358419617830000000}\,\pi\\[0.15cm]
	\approx {} & 0.49999999999998998115\,\pi\,.
\end{align*}
  

\vspace{0.1cm}
\bigskip
\begin{center}
Uwe B\"asel\\[0.15cm]
HTWK Leipzig\\
Fakult\"at Maschinenbau\\
und Energietechnik,\\
04251 Leipzig, Germany\\[0.1cm]
{\small uwe.baesel@htwk-leipzig.de}
\end{center}
\end{document}